\pdfoutput=1
\documentclass[11pt]{amsart}
\usepackage{amsmath,amsthm,amssymb,amsfonts,mathrsfs,mathtools,bm,tikz,amscd}
\usepackage[T1]{fontenc}
\usepackage{lmodern}
\usepackage[utf8]{inputenc}
\usepackage[numbers]{natbib}
\usepackage[backref]{hyperref}
\usepackage{url}
\usetikzlibrary{arrows.meta}
\theoremstyle{definition}
\newtheorem{thm}{Theorem}[subsection]
\newtheorem{dfn}[thm]{Definition}
\newtheorem{lem}[thm]{Lemma}
\newtheorem{cor}[thm]{Corollary}
\newtheorem{prop}[thm]{Proposition}

\newtheorem{rem}[thm]{Remark}

\newcommand{\N}{\mathbb{N}}

\newcommand{\R}{\mathbb{R}}

\newcommand{\per}{\mathrm{Per}}

\newtheorem{introthm}{Theorem}

\title[Existence \& Configuration of Chaotic $\frac{d}{dx}$-Invariant Sets]{Existence and Configuration of Invariant Sets in $C^\infty([a,b])$ on which the Differential Operator Exhibits Devaney's Chaos}

\author{Kazutoyo Iketake}
\address{Department of Mathematics, Faculty of Science, Saitama University, 255 Shimo-Okubo, Sakura-ku, Saitama-shi, Saitama 338-8570, Japan}
\email{iketake2115@gmail.com}

\subjclass[2020]{Primary 47A16; Secondary 37B10, 47B38}
\keywords{Linear chaos, Differential operator, Devaney's chaos, Symbolic dynamics, Hypercyclicity}
\date{\today}

\begin{document}

\begin{abstract}
    In this paper, we investigate the chaotic behavior of the differential operator $\frac{d}{dx}$ on the space of smooth functions $C^\infty([a,b])$ equipped with the $L^p$-norm ($1\le p\le\infty$). We explicitly construct a homeomorphism between a subset of $C^\infty([a,b])$ and the shift space. Moreover, inspired by symbolic dynamics, we demonstrate that invariant sets, on which the differential operator behaves analogously to the shift, are densely configured in $C^\infty([a,b])$. We also prove that the differential operator is chaotic on the entire space $C^\infty([a,b])$ using a similar approach.
\end{abstract}

\maketitle

\section{Introduction}
There has been extensive research on discrete dynamical systems involving linear operators on topological vector spaces. As a typical example, the dynamical behavior of weighted shift operators on sequence spaces is investigated in \cite{MR241956}. Furthermore, pioneering works on Hardy spaces include the study of translation operators by \cite{birkhoff1929demonstration} and differential operators by \cite{MR53231}. In particular, this paper focuses on Devaney's chaos described in \cite{devaney1989introduction}. Since linear operators on finite-dimensional vector spaces cannot be chaotic, linear chaos is an infinite-dimensional phenomenon. In the context of functional analysis, linear chaos on complete spaces, such as Banach spaces and Fréchet spaces, has also been studied. For instance, in \cite{grosse2011linear}, which is a standard text on linear chaos, completeness is assumed in the majority of the arguments. For this reason, standard theoretical tools are not directly applicable to non-complete spaces, and methods tailored to each individual space are required. \par
In this paper, we investigate the chaoticity of the differential operator on the space $C^\infty([a,b])$ of smooth functions on $[a,b]$, equipped with the $L^p$-norm topology. We denote the $L^p$-norm by $\|\cdot\|_p$. Since $(C^\infty([a,b]),\|\cdot\|_p)$ is not complete, applying the general theory is not straightforward. However, it was shown in \cite{MR3918839} that the differential operator is hypercyclic on $(C^\infty([0,1]),\|\cdot\|_\infty)$ (i.e., there exists a dense orbit) using symbolic dynamics and weighted shifts. This implies that the methods well-studied in Hardy spaces and sequence spaces are also effective for incomplete spaces. Building on the idea in \cite{MR3918839}, we analyze the structure and configuration of chaotic invariant subsets of $C^\infty([a,b])$. The main results of this paper consist of the following three theorems. 

\begin{introthm}
    There exists a $\frac{d}{dx}$-invariant subset $A\subset C^\infty([a,b])$ such that $(A,\frac{d}{dx})$ is conjugate to the shift space $(2^{\N},\sigma)$.
\end{introthm}
For a detailed discussion, see \S 4.1 and Theorem~\ref{tekoteko}.

\begin{introthm}
Let \(T:(C^{\infty}([0,b-a]),\|\cdot\|_p) \to (C^{\infty}([a,b]),\|\cdot\|_p)\) be a linear operator defined by $(Tf)(x):=f(x-a)$.  We define $\mathcal{F}=\{F\subset \R\mid 2\le\#F<\infty\}$.
 The family of invariant sets $\{T(E_F)\}_{F\in\mathcal{F}}$ (see Definitions~\ref{def:EF} for the definition) on which the differential operator $\frac{d}{dx}$ is chaotic are densely configured in $C^\infty([a,b])$. 
\end{introthm}

A more general and precise statement is described in Theorem~\ref{kuromichan} and Corollary~\ref{kitexichan}. 

\begin{introthm}
    The discrete dynamical system $(C^\infty([a,b]),\frac{d}{dx})$ is chaotic.
\end{introthm}

This theorem extends the result of \cite[Theorem 3.2]{MR3918839}. The present paper consists of six sections. We prove the main theorems on $C^\infty([0,\gamma])$ in \S 4, and by using results of \S 4,  we prove the main theorems on $C^\infty([a,b])$ in \S 5. In \S 2 and \S 3, we recall some well-known facts and prove several propositions necessary for the proofs of the main theorems.

\section{Preliminaries}
We recall basic concepts in dynamical systems.

\subsection{Chaos and Conjugate}

Let $(X,d)$ be a metric space, and let $f:X\to X$ be a not necessarily continuous map. The pair $(X,f)$ is called a discrete dynamical system. We denote the set of all periodic points of $f$ by $\per(f)$. A set $A\subset X$ is called an $f$-invariant set if it satisfies $f(A) \subseteq A$. When $f$ is clear from the context, we simply call it an invariant set. If $A\subset X$ is an invariant set, then $(A,f|_A)$ is a discrete dynamical system. When the context is clear, we simply write $f$ instead of $f|_A$. \par
Although definitions of chaos vary, we will use the one given by Devaney in \cite{devaney1989introduction}. It was shown in \cite{MR1157223} that the original definition can be simplified to two conditions. Therefore we adopt these two conditions as the definition here.

\begin{dfn}
  Let $A\subset X$ be an invariant set of $f$. The function $f$ is called chaotic (Devaney's chaos) on $(A,d)$ if it satisfies the two conditions below. 
\begin{enumerate}
  \item $\per(f)$ is dense in $(A,d)$. 
  \item $f$ is transitive on $(A,d)$. That is, for any non-empty $(A,d)$-open sets $U,V$, there exists an $n\in\N$ such that $f^n(U)\cap V\ne\emptyset$. 
\end{enumerate}
When no confusion can arise, we simply say that $f$ is chaotic on $A$. 
\end{dfn}

We now recall the definitions and properties of some concepts used in this paper. Our terminology follows one of the standard texts, \cite{grosse2011linear}.

\begin{dfn}
     Let $A\subset X$ be an invariant set of $f$. The function $f$ is called hypercyclic on $(A,d)$ if there exists a point which has a dense orbit. 
\end{dfn}

\begin{dfn}
   Let $f:X_1\to X_1,g:X_2\to X_2$ be maps, and $A_1\subset X_1$ be an invariant set of $f$, $A_2\subset X_2$ be an invariant set of $g$. We say that $f: A_1\to A_1$ and $g: A_2\to  A_2$ are conjugate if there exists a homeomorphism $h:A_1\to A_2$ such that $h\circ f=g\circ h$. The homeomorphism $h$ is called a conjugacy. 
\end{dfn}

\begin{rem}
    $h\circ f=g\circ h$ means that the following diagram is commutative. 
        \[
  \begin{CD}
     { A_1} @>{h}>> { A_2} \\
  @V{f}VV    @V{g}VV \\
     { A_1}   @>{h}>>  { A_2}
  \end{CD}
\]
\end{rem}

\begin{prop}\label{kyouyaku}
    Suppose that $f: A_1\to A_1$ and $g: A_2\to  A_2$ are conjugate. The following four statements hold.
    
    \begin{enumerate}
        \item If $\per(f)$ is dense in $(A_1, d_1)$, then $\per(g)$ is dense in $(A_2, d_2)$.
        \item If $f$ is transitive on $(A_1, d_1)$, then $g$ is transitive on $(A_2, d_2)$. 
        \item If $f$ is hypercyclic on $(A_1, d_1)$, then $g$ is hypercyclic on $(A_2, d_2)$. 
        \item  If $f$ is chaotic on $(A_1, d_1)$, then $g$ is chaotic on $(A_2, d_2)$. 
    \end{enumerate}
\end{prop}

\begin{lem}\label{denseorbit}
    If $f$ is hypercyclic on $A$, then $f$ is transitive on $A$.
\end{lem}

\begin{proof}
    Let $x\in A$ be a point which has a dense orbit and $U,V$ be non-empty open sets of $A$. By density, there exist an $n\in\N$ such that $f^n(x)\in U$ and $m>n$ such that $f^m(x)\in V$. We have $f^m(x)\in f^{m-n}(U)\cap V\ne\emptyset$. 
\end{proof}

\subsection{Symbolic Dynamics}

In this subsection, we recall the standard results in symbolic dynamical systems. Some proofs are omitted since they can be found in many standard textbooks on dynamical systems.

\begin{dfn}
    Let $\Lambda$ be a set defined by $\Lambda:=\{0,1\}^{\N}$.  For $x\in\Lambda$, we denote its $i$-th coordinate by $x_i$. For $x=(x_0,x_1,x_2,\dots),y=(y_0,y_1,y_2,\dots)$, we define maps $d_{\Lambda}:\Lambda\times\Lambda\to\R$ and $\sigma:\Lambda\to\Lambda$ by
    \begin{align*}
      d_{\Lambda}(x,y)&:=\sum_{i=0}^{\infty}\frac{|x_i-y_i|}{2^i}\\
      \sigma(x_0,x_1,x_2,\dots)&:=(x_1,x_2,x_3,\dots). 
    \end{align*}
    $d_\Lambda$ is a metric function on $\Lambda$. Thus, $(\Lambda,d_{\Lambda},\sigma)$ is a discrete dynamical system. $\sigma$ is called the shift map. 
\end{dfn}

\begin{prop}\label{shiftchaos}
    The shift map $\sigma:\Lambda\to\Lambda$ is chaotic on $(\Lambda,d_{\Lambda})$. 
\end{prop}

\begin{proof}
    The proof can be found in Proposition 6.6 of the second edition of \cite{devaney1989introduction}.
\end{proof}
    
\begin{dfn}
  A topological space $C\ne \emptyset$ is called a Cantor set if it satisfies the 
following conditions:
\begin{enumerate}
  \item $C$ is compact. 
  \item $C$ is metrizable. 
  \item $C$ has no isolated points.
  \item $C$ is totally disconnected. 
\end{enumerate}
\end{dfn}

\begin{rem}
    It is well-known that any two Cantor sets are homeomorphic. The proof of this fact can be found in standard texts on descriptive set theory, such as \cite[Theorem 7.4]{MR1321597}. 
\end{rem}

\begin{prop}\label{lmdcantor}
    $(\Lambda,d_{\Lambda})$ is a Cantor set. 
\end{prop}

\section{Auxiliary Propositions and Lemmas}

We provide the proofs of several propositions required for the proof of the main theorems.

\subsection{Metrizability of a Countable Product Topology}

In this subsection, we recall the metrizability theorem for countable product topological spaces. It is well-known that a countable product of metric spaces is metrizable. For example, this fact is presented as an exercise in the standard textbook \cite[\S 10, Exercise 10.]{MR3728284}. The proof of this fact involves the explicit construction of a metric on the product space. However, such a metric can be defined in several ways. Here, we construct the metric in a specific form that is convenient for our subsequent arguments. Furthermore, since it is sufficient to consider the case where the family of metric spaces is bounded, we state the proposition in a specific form.

\begin{dfn}
  Let $(X,d)$ be a metric space. We write $\mathcal{O}(d)$ for the topology on $X$ determined by the metric $d$.
\end{dfn}

\begin{dfn}
  Let $(X,d)$ be a metric space. The diameter of $(X,d)$ is defined by
  \[\mathrm{diam}(X,d)=\mathrm{diam}(X):=\sup_{x,y\in X}d(x,y). \]
  If $\mathrm{diam}(X)<\infty$, we say that $(X,d)$ is bounded. 
\end{dfn}

\begin{prop}\label{pen}
  Let $\{(X_i,d_i)\}_{i=0}^{\infty}$ be a countable family of bounded metric spaces. Assume that $\{\mathrm{diam}(X_i,d_i)\}_{i=0}^{\infty}$ is bounded. Let $(X:=\prod_{i=0}^{\infty}X_i,\mathcal{O})$ be the product space. For any $x,y\in X$, we define a function $d:X\times X\to \R$ by
  \[d(x,y):=\sum_{i=0}^{\infty}\frac{d_i(x_i,y_i)}{2^i},\]
  where $x=(x_0,x_1,\dots), y=(y_0,y_1,\dots)$. Then, $d$ is a metric on $X$ and $\mathcal{O}=\mathcal{O}(d)$. 
\end{prop}

\begin{proof}
    Let $\delta = \sup\{\mathrm{diam}(X_i,d_i)\}_{i=0}^{\infty}$. By assumption, we have $\delta<\infty$. Then 
      \[d(x,y)=\sum_{i=0}^{\infty}\frac{d_i(x_i,y_i)}{2^i}\le\sum_{i=0}^{\infty}\frac{\delta}{2^i}=2\delta<\infty.\]
      Therefore, $d$ is well-defined. It is easy to see that $d$ is a metric. \par
       We will show  $\mathcal{O}=\mathcal{O}(d)$. Let $V\in \mathcal{O}(d)$ and $x\in V$. Since $V$ is open, there exists $\varepsilon>0$ such that $B(x,\varepsilon)\subset V$, where $B(x,\varepsilon)$ is the $\varepsilon$-neighborhood of $x$. Since $\sum\frac{1}{2^i}<\infty$, we can find an $N\in\N$ such that
        \[\sum_{i=N+1}^{\infty}\frac{\delta}{2^i}<\frac{\varepsilon}{2}. \]

Let $\mathrm{pr}_i:(X,\mathcal{O})\to (X_i,\mathcal{O}(d_i))$ be the projection onto the $i$-th coordinate. By the definition of the product topological space, $\mathrm{pr}_i$ is a continuous map. Therefore for each $i=0,1,2,\dots$, there exists an open neighborhood $U_i$ of $x$ in $(X,\mathcal{O})$ such that for any $y\in U_i$, we have
  \[d_i(\mathrm{pr}_i(x), \mathrm{pr}_i(y))=d_i(x_i,y_i)<\frac{2^i\varepsilon}{2(N+1)}.\]
We set $U=U_0\cap U_1\cap\cdots\cap U_N$. Then $U$ is an open neighborhood of $x$ in $(X,\mathcal{O})$. Since for any $y\in U$, 
    \begin{align*}
        d(x,y)&=\sum_{i=0}^{N}\frac{d_i(x_i,y_i)}{2^i}+
        \sum_{i=N+1}^{\infty}\frac{d_i(x_i,y_i)}{2^i}\\
        &<\sum_{i=0}^{N}\frac{1}{2^i}\frac{2^i\varepsilon}{2(N+1)}+
        \sum_{i=N+1}^{\infty}\frac{\delta}{2^i}\\
        &<\frac{\varepsilon}{2}+\frac{\varepsilon}{2}\\
        &=\varepsilon
    \end{align*}
holds and we have $y\in B(x,\varepsilon)$. This implies $U\subset B(x,\varepsilon)\subset  V$ and we get $\mathcal{O}\supset\mathcal{O}(d)$. Conversely, we will show $\mathcal{O}\subset\mathcal{O}(d)$. Let $\mathcal{S}$ be a set defined by
    \[\mathcal{S}=\{\mathrm{pr}_i^{-1}(G)\mid G\in\mathcal{O}(d_i), i=0,1,2,\dots\}.\]
   From the definition of product topology, $\mathcal{S}$ is an open subbase of $(X,\mathcal{O})$. Therefore it suffices to show that $\mathcal{S}\subset \mathcal{O}(d)$.
   Fix an integer $i=0,1,2,\dots$. For any $x,y\in X$, we have 
    \begin{align*}
        d_i(\mathrm{pr}_i(x),\mathrm{pr}_i(y))&=d_i(x_i,y_i)\\
        &=2^i\frac{d_i(x_i,y_i)}{2^i}\\
        &\le 2^i\sum_{j=0}^{\infty}\frac{d_j(x_j,y_j)}{2^j}\\
        &= 2^i d(x,y). 
    \end{align*}
Regarding $\mathrm{pr}_i$ as the map from $(X,d)$ to $(X_i,d_i)$, this inequality shows that the projection $\mathrm{pr}_i$ is a Lipschitz map. Since Lipschitz maps are continuous, the map $\mathrm{pr}_i:(X,\mathcal{O}(d))\to (X_i,\mathcal{O}(d_i))$ is continuous. This means that $\mathrm{pr}_i^{-1}(G)\in\mathcal{O}(d)$ holds for any $G\in \mathcal{O}(d_i)$.  Therefore, it follows that $\mathcal{S}\subset \mathcal{O}(d)$, and the conclusion follows.
\end{proof}

\subsection{Continuity of the Inclusion Map}
The goal of this subsection is to prove Proposition \ref{hiro} and Proposition \ref{eri}. These propositions imply the continuity of the inclusion map and play a key role in the proofs of the main theorems. In the following, we set $\gamma > 0$, $I := [0, \gamma]$ and $1\le p\le\infty$. 

\begin{dfn}
  Let $E\subset C^{\infty}(I)$ be a set defined by
  \[E:=\left\{f\in C^{\infty}(I)\left\lvert f(x)=\sum_{n=0}^{\infty}\frac{a_n}{n!}x^n, a_n\in\{0,1\}\right.\right\}. \]
 For $1\le p \le \infty$, the $L_p$-norm induces a metric on $C^{\infty}(I)$. Thus, $E$  is a metric subspace. In other words, $(E,\rho_p)$ is a metric space with a metric function $\rho_p(f,g):=\|f-g\|_{L^p}$. 
\end{dfn}

\begin{rem}
    It is clear that $f(x)=\sum_{n=0}^{\infty}\frac{a_n}{n!}x^n\in C^{\infty}(I)$ if $a_n\in\{0,1\}$.
\end{rem}

Hereafter, we assume that $f, g \in E$ are given by 
\begin{align*}
    f(x)&=\sum_{n=0}^{\infty}\frac{a_n}{n!}x^n\\ g(x)&=\sum_{n=0}^{\infty}\frac{b_n}{n!}x^n,
\end{align*}
where $a_n, b_n \in \{0,1\}$.

\begin{rem}
  $E$ is not a vector space since it is not closed under addition. Therefore $E$ is not a normed space. 
\end{rem}

The following lemma is known as a corollary of Hölder's inequality. 

\begin{lem}\label{pqineq}
    Let $1\le p < q\le\infty$ and let $\Omega\subset\R$ be a Lebesgue measurable set of finite measure. Then the inequality $\|\cdot\|_p\le |\Omega|^{\frac{1}{p}-\frac{1}{q}}\|\cdot\|_q$ holds, where $|\Omega|$ is the Lebesgue measure of $\Omega$. 
\end{lem}

\begin{proof}
    When $q = \infty$, we have
    \begin{align*}
        \|f\|_p^p=\int_\Omega |f(x)|^pdx\le\left(\sup_{x\in\Omega}|f(x)|\right)^p\int_\Omega dx=|\Omega|\|f\|_\infty^p.
    \end{align*}
   When $q \ne \infty$, applying Hölder's inequality with the exponents $\frac{q}{p}$ and $\frac{q}{q-p}$, we have
\begin{align*}
    \|f\|_p^p=\int_\Omega |f(x)|^pdx\le\left(\int_\Omega |f(x)|^qdx\right)^{\frac{p}{q}}\left(\int_\Omega dx\right)^{\frac{q-p}{q}}=|\Omega|^{\frac{q-p}{q}}\|f\|^p_q.
\end{align*}
Taking the power $\frac{1}{p}$ of both sides gives us the result. 
\end{proof} 

\begin{dfn}
  We write $\frac{d}{dx}:C^{\infty}(I)\to C^{\infty}(I)$ for the differential operator. Since $E$ is an invariant set of $\frac{d}{dx}$, $(E,\rho_p, \frac{d}{dx})$ is a discrete dynamical system.
\end{dfn}

\begin{rem}
 While unbounded operators are often partial maps, all operators in this paper are total maps.
\end{rem}

\begin{dfn}
  For $f,g\in E$, we define a
metric function $d_E$ on $E$ by
  \[d_E(f,g):=\sum_{n=0}^{\infty}\frac{|a_n-b_n|}{(n+1)!}. \]
\end{dfn}

\begin{dfn}
  We define the sequences $\{\eta_k\}_{k=1}^{\infty},\{\zeta_k\}_{k=1}^{\infty},\{\xi_k\}_{k=1}^{\infty}$ as follows:
 \begin{align*}
   \eta_k&:=\sum_{i=k}^{\infty}\frac{1}{i!}\\
      \zeta_k&:=\sum_{i=k}^{\infty}\frac{\gamma^i}{i!}\\
   \xi_k&:=\frac{\gamma^k}{k!}-\sum_{i=k+1}^{\infty}\frac{\gamma^i}{i!}=\frac{\gamma^k}{k!}-\zeta_{k+1}.
 \end{align*}
\end{dfn}

\begin{lem}\label{etaxi}
The following three statements hold.
\begin{enumerate}
    \item \ $\{\eta_k\}_{k=1}^{\infty}, \{\zeta_k\}_{k=1}^{\infty}$ are positive monotonically decreasing sequences. 
    \item $\lim_{k\to\infty}\eta_k=\lim_{k\to\infty}\zeta_k=\lim_{k\to\infty}\xi_k=0$.
    \item There exists $N_{\gamma}\in\N$ such that for any $k\ge N_{\gamma}$, we have $\xi_k>0$ and $\xi_k\ge\xi_{k+1}$. 
\end{enumerate}

\end{lem}

\begin{proof}
  (1) and (2) are clear. We will show (3). Let $N_1\in\N$ with $N_1+2>\gamma$. If $k>N_1$, then
  \begin{align*}
      \sum_{i=k+1}^{\infty}\frac{\gamma^i}{i!}&=\frac{\gamma^{k+1}}{(k+1)!}\left(1+\frac{\gamma}{k+2}+\frac{\gamma^2}{(k+2)(k+3)}+\cdots\right)\\
      &<\frac{\gamma^{k+1}}{(k+1)!}\left(1+\frac{\gamma}{k+2}+\frac{\gamma^2}{(k+2)^2}+\cdots\right)\\
      &=\frac{\gamma^{k+1}}{(k+1)!}\cdot\frac{k+2}{k+2-\gamma}.
  \end{align*}
  Thus, it suffices to show that
\[
\frac{\gamma^k}{k!}>\frac{\gamma^{k+1}}{(k+1)!}\cdot\frac{k+2}{k+2-\gamma}
\]
to show that $\xi_k>0$ for any $k>N_1$. Dividing both sides of this inequality by $\frac{\gamma^k}{k!}$, we have
  \[
1>\frac{\gamma}{k+1}\cdot\frac{k+2}{k+2-\gamma}.
  \]

  Since the right-hand side of this inequality approaches 0 as $k \to \infty$, there exists an $N_2\in\N$ such that for any $k>N_2$, the inequality holds. Therefore, for $k>\max\{N_1,N_2\}$, $\xi_k>0$ holds. \par
  Next, we consider the monotonically decreasing property. For any $k>0$, we have
  \begin{align*}
      \xi_k-\xi_{k+1}&=\frac{\gamma^k}{k!}-\sum_{i=k+1}^{\infty}\frac{\gamma^i}{i!}-\frac{\gamma^{k+1}}{(k+1)!}+\sum_{i=k+2}^{\infty}\frac{\gamma^i}{i!}\\
      &=\frac{\gamma^k}{k!}-2\frac{\gamma^{k+1}}{(k+1)!}\\
      &=\frac{\gamma^k}{k!}\left(1-2\frac{\gamma}{k+1}\right). 
  \end{align*}
 Thus, if we let $N_3\in\N$ with $N_3>2\gamma-1$, then $\xi_k>\xi_{k+1}$ holds for each $k\ge N_3$. Therefore, it suffices to set $N_\gamma=\max\{N_1,N_2,N_3\}$. 
\end{proof}

Henceforth, we fix an $N_\gamma$ obtained from Lemma~\ref{etaxi}.

\begin{lem}\label{apple}
  There exists $M_\gamma\in\N$ such that for any $k\ge M_\gamma$ and $f,g\in E$, if $\rho_1(f,g)<\xi_{k+1}$, then $a_j=b_j$ for each $j=0,1,\dots ,k$. 
\end{lem}

\begin{proof}
For $m\in\N$, we define $\delta_m\in\R$ by

\[\delta_m = \inf\{\rho_1(f,g)\mid f,g\in E,\min\{j\mid a_j\ne b_j\}=m\}.\]

First, we will show that $\delta_m>0$. We define the sequence $\{\alpha_k\}_{k=1}^{\infty}$ by

\[\alpha_k=\frac{1}{k!}-\sum_{i=k+1}^{\infty}\frac{1}{i!}=\frac{1}{k!}-\eta_{k+1}.\]

 For $k\ge1$, 

\begin{align*}
    \eta_{k+1}&=\sum_{i=k+1}^{\infty}\frac{1}{i!}\\
    &=\frac{1}{k!}\sum_{i=1}^{\infty}\frac{1}{(k+1)(k+2)\cdots(k+i)}\\
    &<\frac{1}{k!}\sum_{i=1}^{\infty}\frac{1}{(k+1)^i}\\
    &=\frac{1}{k!}\cdot\frac{1}{k}\le  \frac{1}{k!}.
  \end{align*}
Thus $\alpha_k>0$. Let $f,g\in E$ and let $m=\min\{j\mid a_j\ne b_j\}$. \par
We consider the case $\gamma>1$. We have

\begin{align*}
    \rho_1(f,g)&=\int_0^\gamma \left|\sum_{n=0}^\infty \frac{a_n-b_n}{n!}x^n\right|dx\\
    &\ge \int_0^1 \left|\sum_{n=0}^\infty \frac{a_n-b_n}{n!}x^n\right|dx\\
     &=\int_0^{1}\left|\frac{a_m-b_m}{m!}x^m+ \sum_{n=m+1}^{\infty}\frac{a_n-b_n}{n!}x^n\right|dx\\
    &\ge\int_0^{1}\left(\left|\frac{a_m-b_m}{m!}x^m\right|- \left|\sum_{n=m+1}^{\infty}\frac{a_n-b_n}{n!}x^n\right|\right)dx\\
    &\ge\int_0^{1}\frac{1}{m!}x^mdx-\int_0^{1}\sum_{n=m+1}^{\infty}\frac{|a_n-b_n|}{n!}x^ndx\\
    &=\frac{1}{(m+1)!}-\int_0^{1}\sum_{n=m+1}^{\infty}\frac{|a_n-b_n|}{n!}x^ndx.
\end{align*}
Since $a_n,b_n\in\{0,1\}$, $|a_n-b_n|\le 1$ holds for each $n\in\N$. Noting that the integral and the summation can be exchanged by the Monotone Convergence Theorem, we proceed as follows:
\begin{align*}
    \frac{1}{(m+1)!}-\int_0^{1}\sum_{n=m+1}^{\infty}\frac{|a_n-b_n|}{n!}x^ndx&\ge\frac{1}{(m+1)!}-\int_0^{1}\sum_{n=m+1}^{\infty}\frac{1}{n!}x^ndx\\
   &=\frac{1}{(m+1)!}-\sum_{n=m+1}^{\infty}\frac{1}{n!}\int_0^{1}x^ndx\\
   &=\frac{1}{(m+1)!}-\sum_{n=m+1}^{\infty}\frac{1}{(n+1)!}\\
   &=\alpha_{m+1}.
\end{align*}
Therefore, $\rho_1(f,g)\ge\alpha_{m+1}>0$ holds, and we have $\delta_m>0$. \par
We next consider the case $0<\gamma\le1$. By a similar calculation as before,

\begin{align*}
    \rho_1(f,g)&=\int_0^{\gamma }\left|\frac{a_m-b_m}{m!}x^m+ \sum_{n=m+1}^{\infty}\frac{a_n-b_n}{n!}x^n\right|dx\\
    &\ge\int_0^{\gamma }\left(\left|\frac{a_m-b_m}{m!}x^m\right|- \left|\sum_{n=m+1}^{\infty}\frac{a_n-b_n}{n!}x^n\right|\right)dx\\
    &\ge\int_0^{\gamma }\frac{1}{m!}x^mdx-\int_0^{\gamma }\sum_{n=m+1}^{\infty}\frac{|a_n-b_n|}{n!}x^ndx\\
    &\ge\frac{\gamma^{m+1}}{(m+1)!}-\sum_{n=m+1}^{\infty}\frac{1}{n!}\int_0^{\gamma }x^ndx\\
    &=\frac{\gamma^{m+1}}{(m+1)!}-\sum_{n=m+1}^{\infty}\frac{\gamma^{n+1}}{(n+1)!}\\
    &\ge \frac{\gamma^{m+1}}{(m+1)!}-\sum_{n=m+1}^{\infty}\frac{\gamma^{m+1}}{(n+1)!}\\
    &=\gamma^{m+1}\alpha_{m+1}
\end{align*}
holds. Thus, $\rho_1(f,g)>\gamma^{m+1}\alpha_{m+1}>0$ holds, and we have $\delta_m>0$. \par
We define $\delta\in \R$ by
\[\delta=\min_{m\in\{0,1,\dots, N_\gamma-1\}}\delta_m.\]
 As $\delta_m>0$, it follows that $\delta>0$. By Lemma~\ref{etaxi}, we can choose $N_0\in\N$ such that for any $n>N_0$, we have $\xi_{n+1}\le\delta$. Let $M_\gamma=\max\{N_0+1,N_\gamma\}$. We fix $k\ge M_\gamma$ and $f,g\in E$. \par
 Henceforth we will show the contrapositive statement:
\[\text{If there exists }j\in\{0,1,\dots ,k\}\text{ such that }a_j\ne b_j\text{, then }\rho_1(f,g)\ge\xi_{k+1}. \]
Let $l\in\{0,1,\dots k\}$ be the smallest integer such that $a_l\ne b_l$. When $l<M_\gamma$, it follows from the definitions of $\delta$ and $N_0$ that
\begin{align*}
    \rho_1(f,g)&\ge\delta_l\ge\delta\ge\xi_{k+1}.
\end{align*}
  When $l\ge  M_\gamma$, we have
  \begin{align*}
    \rho_1(f,g)
    &=\int_0^{\gamma}\left|\frac{a_l-b_l}{l!}x^l+ \sum_{n=l+1}^{\infty}\frac{a_n-b_n}{n!}x^n\right|dx\\
    &\ge\int_0^{\gamma}\left(\left|\frac{a_l-b_l}{l!}x^l\right|- \left|\sum_{n=l+1}^{\infty}\frac{a_n-b_n}{n!}x^n\right|\right)dx\\
     &\ge\frac{\gamma^{l+1}}{(l+1)!}-\int_0^{\gamma}\sum_{n=l+1}^{\infty}\frac{1}{n!}x^ndx\\
   &=\frac{\gamma^{l+1}}{(l+1)!}-\sum_{n=l+1}^{\infty}\frac{1}{n!}\int_0^{\gamma}x^ndx\\
   &=\frac{\gamma^{l+1}}{(l+1)!}-\sum_{n=l+1}^{\infty}\frac{\gamma^{n+1}}{(n+1)!}\\
   &=\xi_{l+1}.
  \end{align*}
Since $l\ge  M_\gamma\ge N_\gamma$, from Lemma~\ref{etaxi}, $\rho_1(f,g)\ge\xi_{l+1}\ge\xi_{k+1}$ holds. 
\end{proof}

Henceforth, we fix an $M_\gamma$ obtained from Lemma~\ref{apple}.

\begin{lem}\label{banana}
   Let $k\ge1$ be an integer. For any $f,g\in E$, if $a_j=b_j$ for $j=0,1,\dots,k$, then $d_E(f,g)\le\eta_{k+2}$. 
\end{lem}

\begin{proof}
Suppose that $a_j=b_j$ for each $j=0,1,\dots,k$. We have 
    \begin{align*}
        d_E(f,g)&=\sum_{n=0}^{\infty}\frac{|a_n-b_n|}{(n+1)!}=\sum_{n=k+1}^{\infty}\frac{|a_n-b_n|}{(n+1)!}\le\sum_{n=k+1}^{\infty}\frac{1}{(n+1)!}=\eta_{k+2},  
    \end{align*}
  as required.
\end{proof}

\begin{prop}\label{hiro}
   For any $\varepsilon>0$ and $f\in E$, there exists a $\delta>0$ such that for any $g\in E$ with $\rho_1(f,g)<\delta$, we have $d_E(f,g)<\varepsilon$.
\end{prop}

\begin{proof}
    Let $\varepsilon>0$ and let $f\in E$. By Lemma~\ref{etaxi}, we can find an $N\ge M_\gamma$ such that $\eta_N<\varepsilon$. 
    Let $\delta=\xi_{N+1}$, and fix $g\in E$ with $\rho_1(f,g)<\delta$. From Lemma~\ref{apple}, for each $j=0,1,\dots,N$, we have $a_j=b_j$. Lemma~\ref{banana} and Lemma~\ref{etaxi} give us 
    \[d_E(f,g)\le\eta_{N+2}\le\eta_{N}<\varepsilon, \]
        and the conclusion follows.
\end{proof}

\begin{lem}\label{orange}
 Let $k\ge1$ be an integer. For any $f,g\in E$, if $d_E(f,g)<\frac{1}{(k+1)!}$, then $a_j=b_j$ for $j=0,1,\dots ,k$. 
\end{lem}

\begin{proof}
  We prove the contrapositive. Let $l\in\{0,1,\dots k\}$ be the smallest integer such that $a_l\ne b_l$. Since $a_l,b_l\in\{0,1\}$, we have $|a_l-b_l|=1$. Therefore
    \[d_E(f,g)=\sum_{n=0}^{\infty}\frac{|a_n-b_n|}{(n+1)!}=\sum_{n=l}^{\infty}\frac{|a_n-b_n|}{(n+1)!}\ge\frac{1}{(l+1)!}\ge\frac{1}{(k+1)!}, \]
    and the conclusion follows.
\end{proof}

\begin{lem}\label{melon}
   Let $k\ge1$ be an integer. For any $f,g\in E$, if $a_j=b_j$ for any $j=0,1,\dots,k$, then $\rho_{\infty}(f,g)\le\zeta_{k+1}$. 
\end{lem}

\begin{proof}
    Suppose that $a_j=b_j$ for each $j=0,1,\dots,k$. We have 
    \begin{align*}
        \rho_{\infty}(f,g)&=\sup_{x\in I}\left| \sum_{n=0}^{\infty}\frac{a_n-b_n}{n!}x^n \right|\\
        &=\sup_{x\in I}\left| \sum_{n=k+1}^{\infty}\frac{a_n-b_n}{n!}x^n \right|\\
        &\le\sup_{x\in I}\sum_{n=k+1}^{\infty}\frac{|a_n-b_n|}{n!}x^n \\
        &\le\sup_{x\in I}\sum_{n=k+1}^{\infty}\frac{1}{n!}x^n \\
         &=\sum_{n=k+1}^{\infty}\frac{\gamma^n}{n!}\\
         &=\zeta_{k+1}, 
    \end{align*}
  as required.
\end{proof}

\begin{prop}\label{eri}
For any $\varepsilon>0$ and $f\in E$, there exists a $\delta>0$ such that for any $g\in E$ with $d_E(f,g)<\delta$, we have $\rho_{\infty}(f,g)<\varepsilon$.
\end{prop}

\begin{proof}
    Let $\varepsilon>0$ and let $f\in E$. By Lemma~\ref{etaxi}, we can find an $N\ge1$ such that $\zeta_{N}<\varepsilon$. 
    Let $\delta=\frac{1}{(N+1)!}$, and fix $g\in E$ with $d_E(f,g)<\delta$. From Lemma~\ref{orange}, for each $j=0,1,\dots,N$, we have $a_j=b_j$. Lemma~\ref{etaxi} and  Lemma~\ref{melon} give us 
    \[\rho_\infty(f,g)\le\zeta_{N+1}\le\zeta_{N}<\varepsilon, \]
        and the conclusion follows.
\end{proof}

\section{Proof of the Main Theorems on $[0,\gamma]$}

In this section, we prove the main theorems on $C^\infty{([0,\gamma])}$. The proof on $C^\infty{([a,b])}$ will be given in the next section. We recall that $I=[0,\gamma]$.

\subsection{Existence of an invariant set in $C^\infty{([0,\gamma])}$ conjugate to a shift space}

The set $E$, defined in the previous section, has a natural correspondence with the set $\Lambda$, which is defined in \S 2.2. Thus, we aim to show that $(E, \frac{d}{dx})$ is conjugate to $(\Lambda, \sigma)$. While the commutativity condition in the definition of conjugacy is straightforward to verify, it is non-trivial to show that the correspondence is a homeomorphism. We prepare several lemmas to demonstrate that these two sets are equipped with equivalent metric structures.

\begin{dfn}
    For each $i=0,1,2,\dots$, let $D_i:=\{0,1\}$. For any $x,y\in D_i$, we define metric functions $d_i^{(1)},d_i^{(2)}:D_i\times D_i\to \R$ as follows:
    \begin{align*}
  d_i^{(1)}(x,y)&:=|x-y|=\left\{ \,
    \begin{aligned}
    & 1\ &(x\ne y)  \\
    & 0\ &(x= y)  \\
    \end{aligned}
\right.\\
  d_i^{(2)}(x,y)&:=\frac{2^i}{(i+1)!}|x-y|
=\left\{ \,
    \begin{aligned}
    & \frac{2^i}{(i+1)!}\ &(x\ne y)  \\
    & 0\ &(x= y)  \\
    \end{aligned}
\right. .
    \end{align*}
\end{dfn}

\begin{rem}
       By the definition of the Cartesian product, $\prod_{i=0}^{\infty}D_i=\Lambda$. 
\end{rem}

\begin{dfn}
    For $m=1,2$, we write $\mathcal{O}_m$ for the product topology of the family $(D_i,d_i^{(m)})_{i=0}^{\infty}$. 
\end{dfn}

\begin{lem}\label{lmdhomeo}
   The topological spaces $(\Lambda,\mathcal{O}(d_{\Lambda})),(\Lambda,\mathcal{O}_1)$ and $(\Lambda,\mathcal{O}_2)$ are homeomorphic to each other. Even stronger, $\mathcal{O}(d_{\Lambda})=\mathcal{O}_1=\mathcal{O}_2$. 
\end{lem}

\begin{proof}
    First, we will show that $\mathcal{O}(d_{\Lambda})=\mathcal{O}_1$. For any $x,y\in\Lambda$, we define $d:\Lambda\times\Lambda\to\R$ by
       \[d(x,y):=\sum_{i=0}^{\infty}\frac{d_i^{(1)}(x,y)}{2^i}. \]
    By Proposition~\ref{pen}, we have $\mathcal{O}_1=\mathcal{O}(d)$. From the definition of $d_i^{(1)}$, 
\begin{align*}
        d(x,y)=\sum_{i=0}^{\infty}\frac{d_i^{(1)}(x,y)}{2^i}=\sum_{i=0}^{\infty}\frac{|x_i-y_i|}{2^i}=d_{\Lambda}(x,y).
\end{align*}
This means $d=d_{\Lambda}$, therefore  $\mathcal{O}(d_{\Lambda})=\mathcal{O}_1$ holds.\par 
Next, we will show that $\mathcal{O}_1=\mathcal{O}_2$. For each $i=0,1,2,\dots$, both $d_i^{(1)}$ and $d_i^{(2)}$ induce the discrete topology on $D_i$. Proposition~\ref{pen} implies $\mathcal{O}_1=\mathcal{O}_2$.  
\end{proof}

\begin{lem}\label{rhohomeo}
    The topological spaces $(E,\mathcal{O}(\rho_p))$ and $(E,\mathcal{O}(d_E))$ are homeomorphic. Even stronger, $\mathcal{O}(\rho_p)=\mathcal{O}(d_E)$. 
\end{lem}

\begin{proof}
We denote the open $\varepsilon$-neighborhoods of $f$ in $(E,\rho_p)$ and $(E,d_E)$ by $B_p(f,\varepsilon)$ and $B_E(f,\varepsilon)$, respectively. \par
Let $U\in\mathcal{O}(\rho_p)$ and let $f\in U$. Since $U$ is open, there exists $\varepsilon>0$ such that $B_p(f,\varepsilon)\subset U$. By Proposition~\ref{eri}, we can find a $\delta>0$ such that for any $g\in E$ with $d_E(f,g)<\delta$, we have $\rho_{\infty}(f,g)<\gamma^{-\frac{1}{p}}\varepsilon$. In other words, 
\[B_E(f,\delta)\subset B_{\infty}(f,\gamma^{-\frac{1}{p}}\varepsilon). \]
From Lemma~\ref{pqineq}, we have
\[B_{\infty}(f,\gamma^{-\frac{1}{p}}\varepsilon)\subset B_p(f,\varepsilon)\]
Therefore $B_E(f,\delta)\subset B_p(f,\varepsilon)\subset U$, and this means  $U\in \mathcal{O}(d_E)$. Thus we have $\mathcal{O}(\rho_p)\subset\mathcal{O}(d_E)$. 
The reverse inclusion $\mathcal{O}(\rho_p)\supset\mathcal{O}(d_E)$ similarly follows from Proposition~\ref{hiro} and Lemma~\ref{pqineq}. 
\end{proof}

The following theorem is the first main theorem on $C^\infty(I)$

\begin{thm}\label{potato}
   $(\Lambda,d_{\Lambda},\sigma)$ and $(E,\rho_p, \frac{d}{dx})$ are conjugate.
\end{thm}

\begin{proof}
   From Lemma~\ref{lmdhomeo} and Lemma~\ref{rhohomeo}, the identity maps
    \[\mathrm{id}_{\Lambda}:(\Lambda,\mathcal{O}(d_{\Lambda}))\to(\Lambda,\mathcal{O}_2)\]
   and 
   \[\mathrm{id}_{E}:(E,\mathcal{O}(d_E))\to(E,\mathcal{O}(\rho_p))\]
   are homeomorphisms. \par
   For $a=(a_0,a_1,a_2,\dots)\in \Lambda$, we define a map $\iota:(\Lambda,\mathcal{O}_2)\to (E,\mathcal{O}(d_E))$ by
       \[\iota(a)(x):=\sum_{n=0}^{\infty}\frac{a_n}{n!}x^n\quad(x\in I). \]
  It is clear that $\iota$ is a bijection. Let the metric $d$ on $\Lambda$ be defined by
    \[d(x,y)=\sum_{i=0}^{\infty}\frac{d_i^{(2)}(x,y)}{2^i},\]
Proposition~\ref{pen} gives us $\mathcal{O}_2=\mathcal{O}(d)$. Thus we can regard $\iota$ as a continuous map from $(\Lambda,d)$ to $(E,d_E)$. 
 For any $a=(a_0,a_1,a_2,\dots),b=(b_0,b_1,b_2,\dots)\in \Lambda$, we have
    \begin{align*}
        d_E(\iota(a),\iota(b))&=\sum_{n=0}^{\infty}\frac{|a_n-b_n|}{(n+1)!}\\
        &=\sum_{i=0}^{\infty}\frac{1}{2^i}\frac{2^i}{(i+1)!}|a_i-b_i|\\
        &=\sum_{i=0}^{\infty}\frac{d_i^{(2)}(a_i,b_i)}{2^i}\\
        &=d(a,b). 
    \end{align*}
    This calculation implies that $\iota:(\Lambda,d)\to (E,d_E)$ is an isometry. 
    Since a surjective isometry is a homeomorphism, $\iota:(\Lambda,\mathcal{O}_2)\to (E,\mathcal{O}(d_E))$ is also a homeomorphism. Therefore $\mathrm{id}_E\circ\iota \circ \mathrm{id}_{\Lambda}$ is a homeomorphism. \par
    For any $a=(a_0,a_1,a_2,\dots)\in \Lambda$ and $x\in I$, we have
    \begin{align*}
        \iota(\sigma(a))(x)&=\sum_{n=0}^{\infty}\frac{a_{n+1}}{n!}x^n\\
        &=\frac{d}{dx}\left(\sum_{n=0}^{\infty}\frac{a_{n}}{n!}x^n\right)\\
        &=\frac{d}{dx}\iota(a)(x). 
    \end{align*}
    This means that $\iota\circ\sigma=\frac{d}{dx}\circ\iota$. It is clear that $\mathrm{id}_{\Lambda}\circ \sigma =\sigma\circ \mathrm{id}_{\Lambda}$ and $\mathrm{id}_{E}\circ\frac{d}{dx}=\frac{d}{dx}\circ \mathrm{id}_{E}$.\par
    Therefore, the following diagram is commutative, and $\mathrm{id}_E\circ\iota \circ \mathrm{id}_{\Lambda}$ is a conjugacy. 
  
            \[
  \begin{CD}
     { (\Lambda,\mathcal{O}(d_{\Lambda}))} @>{\mathrm{id}_{\Lambda}}>> { (\Lambda,\mathcal{O}_2)} @>{\iota}>> { (E,\mathcal{O}(d_E))}@>{\mathrm{id}_{E}}>> { (E,\mathcal{O}(\rho_p))} \\
  @V{\sigma}VV    @V{\sigma}VV  @V{\frac{d}{dx}}VV  @V{\frac{d}{dx}}VV\\
     { (\Lambda,\mathcal{O}(d_{\Lambda}))}   @>{\mathrm{id}_{\Lambda}}>>  { (\Lambda,\mathcal{O}_2)} @>{\iota}>>  { (E,\mathcal{O}(d_E))} @>{\mathrm{id}_{E}}>>  { (E,\mathcal{O}(\rho_p))}
  \end{CD}
\]

\end{proof}

\begin{cor}\label{togetoge}
The following two statements hold.
\begin{enumerate}
    \item   The differential operator $\frac{d}{dx}$ is chaotic on $(E,\rho_p)$.
    \item   $E$ is a Cantor set. 
\end{enumerate}
\end{cor}

\begin{proof}
    Proposition~\ref{kyouyaku}, Proposition~\ref{shiftchaos} and Theorem~\ref{potato} imply (1). (2) follows from Theorem~\ref{potato}, Proposition~\ref{lmdcantor} and the fact that all conditions in the definition of a Cantor set are topological properties.
\end{proof}

\subsection{Configuration of Chaotic Invariant Sets}

In Corollary~\ref{togetoge}, we have shown that $E$ is a chaotic invariant set. In this subsection, we consider sets that are similar to $E$ and their configuration. \\

Let $F=\{c_0,c_1,\dots,c_m\} \subset \R$ be a finite set with $m \ge 1$. The diameter of $F$ is given by 
\begin{align*}
    \mathrm{diam}(F):=\max_{i,j\in\{0,1,\dots,m\}}|c_i-c_j|. 
\end{align*}
In the following discussion, we let $C^\infty(I)$ denote the set of $C^\infty$ real-valued functions, $C^\infty(I)=C^\infty(I,\R)$. However, our arguments also hold if $C^\infty(I)$ is taken as the set of $C^\infty$ complex-valued functions, $C^\infty(I)=C^\infty(I,\mathbb{C})$. In the latter case, $F \subset \mathbb{C}$.

\begin{dfn}\label{def:EF}
      Let $E_{F}\subset C^{\infty}(I)$ be the set defined by
  \[E_{F}:=\left\{f\in C^{\infty}(I)\left\lvert f(x)=\sum_{n=0}^{\infty}\frac{a_n}{n!}x^n, a_n\in F\right.\right\}.\]
   For $1\le p \le \infty$, $(E_F,\rho_p)$ is a metric space. $E_F$ is an invariant set of $\frac{d}{dx}$. Thus, $(E_F,\rho_p, \frac{d}{dx})$ is a discrete dynamical system. 
\end{dfn}

\begin{prop}\label{renzokutainoudo}
    $E_F$ has the cardinality of the continuum; that is, $|E_F|=|\R|$, where $|\cdot|$ denotes the cardinality of a set.
\end{prop}

\begin{proof}
    Recall $F=\{c_0,c_1,\dots,c_m\}$. We define a map 
    \[\varphi:E_F\to\{0,1,2,\dots,m\}^\N=(m+1)^\N\]
    by
    \begin{align*}
        \mathrm{pr}_i\circ\varphi\left( \sum_{n=0}^{\infty}\frac{a_n}{n!}x^n\right)=j\quad \iff\quad a_i=c_j
    \end{align*}
    for $i=0,1,2,\dots$, where $\mathrm{pr}_i$ is the $i$-th projection map. It is clear that $\varphi$ is a bijection. Thus, 
    \begin{align*}
        |E_F|=|(m+1)^\N|=|\R|
    \end{align*}
    holds.
\end{proof}

Hereafter, we assume that $f, g \in E_F$ are given by 
\begin{align*}
    f(x)=\sum_{n=0}^{\infty}\frac{a_n}{n!}x^n\\ g(x)=\sum_{n=0}^{\infty}\frac{b_n}{n!}x^n,
\end{align*}
where $a_n, b_n \in F$.

\begin{lem}\label{piano}
    Let $k\ge1$ be an integer. For any $f,g\in E_F$, if $a_j=b_j$ for $j=0,1,\dots,k$, then $\rho_{\infty}(f,g)\le \mathrm{diam}(F)\zeta_{k+1}$.
\end{lem}

\begin{proof}
    This can be shown similarly to Lemma~\ref{melon} by using   $|a_n-b_n|\le \mathrm{diam}(F)$.
\end{proof}

\begin{thm}\label{kao}
The differential operator $\frac{d}{dx}$ is chaotic on $(E_F, \rho_p)$. 
\end{thm}

\begin{proof}
    First, we will show that $\per(\frac{d}{dx})$ is dense in $(E_F, \rho_p)$. Let $f\in E_F$ and $\varepsilon>0$. By Lemma~\ref{etaxi}, we can find an $N\ge1$ such that $\gamma^{\frac{1}{p}}\mathrm{diam}(F)\zeta_N<\varepsilon$. Let $f(x)=\sum_{n=0}^{\infty}\frac{a_n}{n!}x^n$. 
For $j=0,1,2,\dots$, we define $b_j$ by
\begin{align*}
     b_j&=\left\{ \,
    \begin{aligned}
    & a_j\ &(j=0,1,\dots,N)  \\
    & b_{j-(N+1)}\ &(j\ge N+1)  \\
    \end{aligned}
\right. 
\end{align*}
  and $g(x)=\sum_{n=0}^{\infty}\frac{b_n}{n!}x^n$. It is clear that $g$ is an $(N+1)$-periodic point of $\frac{d}{dx}$. From Lemma~\ref{piano}, Lemma~\ref{pqineq} and Lemma~\ref{etaxi}, we have
  \[\rho_p(f,g)\le\gamma^{\frac{1}{p}}\rho_{\infty}(f,g)\le \gamma^{\frac{1}{p}}\mathrm{diam}(F)\zeta_{N+1}\le \gamma^{\frac{1}{p}}\mathrm{diam}(F)\zeta_{N}<\varepsilon\]
  which implies the density of $\per(\frac{d}{dx})$. \par
    Next, we will show transitivity. By Lemma~\ref{denseorbit}, it is sufficient to show the existence of a point in $E_F$ with a dense orbit. 
Let $L:=\{(s_1,s_2,\dots,s_k) \mid k\ge 1, s_n\in F\}$ be the set of all finite sequences of elements from $F$. Since $F$ is finite, $L$ is countable. 
Let $b=(b_0,b_1,b_2, \dots)\in F^{\N_{\ge0}}$ be the sequence of all elements of $L$. For example, if $F=\{0,1,2\}$, then
  \begin{align*}
      b=(0\ 1\ 2\mid 00\ 01\ 02\ 10\ 11\ 12\ 20\ 21\ 22\mid 000\ 001\ 002\ 100\ 101\ \dots). 
  \end{align*}
Let $g\in E_F$ be a function defined by $g(x)=\sum_{n=0}^{\infty}\frac{b_n}{n!}x^n$. We will show that the orbit of $g$ is dense. Let $f\in E_F$, and $\varepsilon>0$. From Lemma~\ref{etaxi}, we can find an $N\ge1$ such that $\gamma^{\frac{1}{p}}\mathrm{diam}(F)\zeta_N<\varepsilon$.
Let $f(x)=\sum_{n=0}^{\infty}\frac{a_n}{n!}x^n$. By the definition of $g$, there exists $l\in\N$ such that
\[ b_{l}=a_0,b_{l+1}=a_1,\dots ,b_{l+N}=a_N.\]
Then, 
  \begin{align*}
      \frac{d^l}{dx^l}g(x)=\sum_{n=0}^{\infty}\frac{b_{n+l}}{n!}x^n. 
  \end{align*}
This and Lemma~\ref{piano}, Lemma~\ref{pqineq} and Lemma~\ref{etaxi} imply
   \[\rho_p\left(f,\frac{d^l}{dx^l}g\right)\le\gamma^{\frac{1}{p}}\rho_{\infty}\left(f,\frac{d^l}{dx^l}g\right)\le \gamma^{\frac{1}{p}}\mathrm{diam}(F)\zeta_{N+1}\le \gamma^{\frac{1}{p}}\mathrm{diam}(F)\zeta_{N}<\varepsilon.\]
  Therefore the orbit of $g$ is dense in $(E_F,\rho_p)$.
\end{proof}

\begin{rem}
   Presumably, analogous to Theorem \ref{potato}, we can show that $(E_F, \frac{d}{dx})$ is conjugate to a shift space. However, to avoid unnecessary complexity and to extract the essential arguments, we do not provide this proof in this paper.
\end{rem}

The following theorem is well known as the Weierstrass Approximation Theorem.

\begin{prop}\label{natu}
 Let $f:\Omega\to\mathbb{K}$ be a continuous function, where $\mathbb{K}=\R$ or $\mathbb{C}$ and  $\Omega$ is an interval or $\R$. Then, for any $\varepsilon>0$ there exists a polynomial function $P_{\varepsilon}$ such that $\rho_\infty(f,P_{\varepsilon})<\varepsilon$. 
\end{prop}

\begin{dfn}
    Let $P$ be a polynomial function given by
    \[P(x)=\frac{p_n}{n!}x^n+\frac{p_{n-1}}{(n-1)!}x^{n-1}+\cdots+\frac{p_1}{1!}x+\frac{p_0}{0!}.\]
  We define the set $F(P)$ by
    \[F(P)=\{0,p_0,p_1,\dots ,p_n\}.\]
\end{dfn}

We will slightly modify the Weierstrass Approximation Theorem into a form that can be used in the proof of the main theorems.

\begin{prop}\label{waieru}
Let $f\in C(I)$. For any $\varepsilon>0$, there exists a polynomial function $P_\varepsilon$ satisfying the following two conditions.
\begin{enumerate}
    \item $\rho_\infty(f,P_{\varepsilon})<\varepsilon$. 
    \item $\#F(P_\varepsilon)\ge2$.
\end{enumerate}
\end{prop}

\begin{proof}
 Let $f\in C(I)$ and let $\varepsilon>0$. By Proposition~\ref{natu} we can find a polynomial function $P_{1/2\varepsilon}$ such that $\rho_\infty(f,P_{1/2\varepsilon})<\frac{1}{2}\varepsilon$. If $\#F(P_{1/2\varepsilon})>1$, the result holds trivially. We consider the case $\#F(P_{1/2\varepsilon})=1$. Since $F(P_{1/2\varepsilon})=\{0\}$ in this case, we have $P_{1/2\varepsilon}=0$. We choose an $\alpha\in\R$ satisfying $0<\alpha<\frac{1}{2}\varepsilon$. Then, 
 \begin{align*}
\rho_\infty(P_{1/2\varepsilon},P_{1/2\varepsilon}+\alpha)&=\|\alpha\|_\infty=\sup_{x\in I}|\alpha|=\alpha<\frac{1}{2}\varepsilon
 \end{align*}
     holds. By the triangle inequality, we have
     \[\rho_\infty(f,P_{1/2\varepsilon}+\alpha)\le\rho_\infty(f,P_{1/2\varepsilon})+\rho_\infty(P_{1/2\varepsilon},P_{1/2\varepsilon}+\alpha)<\frac{1}{2}\varepsilon+\frac{1}{2}\varepsilon=\varepsilon.\]
Thus, if we set $P_\varepsilon=P_{1/2\varepsilon}+\alpha$, then $F(P_\varepsilon)=\{0,\alpha\}$, as required.    
\end{proof}

The following is the second main theorem. It means that chaotic invariant sets are densely configured, where we exclude trivially chaotic finite invariant sets.

\begin{thm}\label{syut}
    For any $f\in C^{\infty}(I)$ and $\varepsilon>0$, there exists a finite set  $F\subset\R$ with $2\le \#F <\infty$ such that $\rho_p(f,E_F)<\varepsilon$. 
\end{thm}

\begin{proof}
    Let $f\in C^{\infty}(I)$ and let $\varepsilon>0$. By Lemma~\ref{waieru}, we can find a polynomial function $P$ such that we have $\rho_\infty(f,P)<\gamma^{-\frac{1}{p}}\varepsilon$ and $\#F(P)\ge 2$. From Theorem~\ref{kao} with $F=F(P)$, $\frac{d}{dx}$ is chaotic on $(E_F,\rho_p)$. By Lemma~\ref{pqineq} and $P\in E_F$, we have
    \[\rho_p(f,E_F)=\inf_{g\in E_F}\rho_p(f,g)\le\inf_{g\in E_F}\gamma^{\frac{1}{p}}\rho_{\infty}(f,g)\le\gamma^{\frac{1}{p}}\rho_{\infty}(f,P)<\gamma^{\frac{1}{p}}\gamma^{-\frac{1}{p}}\varepsilon=\varepsilon,\]
    as required.
    
\end{proof}

\begin{cor}\label{shimauma}
     For any $f\in C^{\infty}(I)$, there exists a filtration $\{E_{F(n)}\}_{n=1}^\infty$ approximating $f$. That is, there exists a family $\{F(n)\}_{n=1}^{\infty}\subset\R$ satisfying the following conditions.
    \begin{enumerate}
        \item $2\le\#F(n)<\infty$ for any $n\ge 1$
        \item $E_{F(n)}\subset E_{F(n+1)}$ for any $n\ge1$.
        \item $f\in\overline{\bigcup_{n=1}^\infty E_{F(n)}}$, where the closure is with respect to $\rho_p$.
    \end{enumerate}
\end{cor}

\begin{proof}
    We fix $f\in C^{\infty}(I)$. Let $n\ge1$ be an integer. By Theorem~\ref{syut}, there exists $F'(n)\subset\R$ such that $2\le\#F'(n)<\infty$ and $\rho_p(f,E_{F'(n)})<\frac{1}{n}$. We define a set $F(n)$ by
        \[F(n)=\bigcup_{k=1}^{n}F'(k).\]
    We consider the family $\{E_{F(n)}\}_{n=1}^{\infty}$. It is clear that $2\le\#F(n)<\infty$ and $E_{F(n)}\subset E_{F(n+1)}$ for any $n\ge1$. Let $\varepsilon>0$. Then, we can find $N\in\N$ such that $\frac{1}{N}<\varepsilon$. Since $F'(N)\subset F(N)$, we have $E_{F'(N)}\subset E_{F(N)}$. Thus,
        \begin{align*}
            \rho_p(f,E_{F(N)})\le\rho_p(f,E_{F'(N)})<\frac{1}{N}<\varepsilon
        \end{align*}
    holds. 
\end{proof}

\subsection{Chaos of $(C^\infty{(I)}, \frac{d}{dx})$}

While the previous sections focused on chaotic proper subsets of $C^\infty{(I)}=C^\infty{([0,\gamma])}$, this subsection investigates the chaoticity of $C^\infty{(I)}$ itself. It was shown in \cite{MR3918839} that $(C^\infty{([0,1])}, \rho_\infty, \frac{d}{dx})$ is hypercyclic. By a similar method, it can be proven that $(C^\infty{(I)}, \rho_\infty,\frac{d}{dx})$ is also hypercyclic. The following proposition is obtained from Theorem 3.2 of \cite{MR3918839}.

\begin{prop}\label{mio}
    The differential operator $\frac{d}{dx}$ is hypercyclic on the metric space $(C^\infty(I),\rho_\infty)$.
\end{prop}

Although the result in \cite{MR3918839} was for the $L^\infty$-norm, it can be easily extended to the $L^p$-norm.

\begin{prop}\label{haru}
   The differential operator $\frac{d}{dx}$ is hypercyclic on the metric space $(C^\infty(I),\rho_p)$.
\end{prop}

\begin{proof}
    By Proposition~\ref{mio}, there exists $f\in C^\infty(I)$ such that $f$ has a dense orbit in $(C^\infty(I),\rho_\infty)$. Let $g\in C^\infty(I)$ and let $\varepsilon>0$. By the definition of $f$, there exists an $N\in\N$ such that 
    \[\rho_\infty\left(\frac{d^N}{dx^N}f,g\right)<\gamma^{-\frac{1}{p}}\varepsilon.\]
    From Lemma~\ref{pqineq}, we have
    \[\rho_p\left(\frac{d^N}{dx^N}f,g\right)\le\gamma^\frac{1}{p}\rho_\infty\left(\frac{d^N}{dx^N}f,g\right)<\gamma^\frac{1}{p}\gamma^{-\frac{1}{p}}\varepsilon=\varepsilon,\]
    and the conclusion follows. 
\end{proof}

The transitivity in the definition of chaos follows from Lemma~\ref{denseorbit}. Next, we prove the density of periodic points.

\begin{prop}\label{huyu}
    $\per(\frac{d}{dx})$ is dense in $(C^\infty(I), \rho_p)$.
\end{prop}

\begin{proof}
    Let $f\in C^\infty(I)$ and let $\varepsilon>0$. By Lemma~\ref{waieru}, we can find a polynomial function with
        \[\rho_\infty(f,P)<\frac{1}{2}\gamma^{-\frac{1}{p}}\varepsilon,\quad \#F(P)\ge 2.\]
    We write $P$ as
        \[P(x)=\frac{p_n}{n!}x^n+\frac{p_{n-1}}{(n-1)!}x^{n-1}+\cdots+\frac{p_1}{1!}x+\frac{p_0}{0!},\]
    where $n=\deg(P)$. By Lemma~\ref{etaxi}, there exists an $N\ge 1$ such that
        \[\zeta_N<\frac{\gamma^{-\frac{1}{p}}}{2\ \mathrm{diam}(F(P))}\varepsilon.\]
    We set $L=\max\{0,N-n\},L'=\max\{n,N\}$. For each $j=0,1,2,\dots$, we define $b_j$ by
\begin{align*}
     b_j&:=\left\{ \,
    \begin{aligned}
    & p_j\ &(0\le j \le n) \\
    &  0   &    (j\in \{k\in\N\mid n< k \le n+L\})  \\
    & b_{j-L'-1}\ &(j> L')  \\
    \end{aligned}
\right.
\end{align*}
and we define $g\in C^\infty(I)$ by $g(x)=\sum_{k=0}^\infty\frac{b_k}{k!}x^k$.  Since $g,P\in E_{F(P)}$, from Lemma~\ref{pqineq}, Lemma~\ref{etaxi} and Lemma~\ref{piano}, we have

\begin{align*}
    \rho_p(P,g)&\le\gamma^\frac{1}{p}\rho_{\infty}(P,g)\\
    &\le \gamma^\frac{1}{p}\mathrm{diam}(F(P))\zeta_{L'+1}\\
    &\le \gamma^\frac{1}{p}\mathrm{diam}(F(P))\zeta_{N}\\
    &< \gamma^\frac{1}{p}\mathrm{diam}(F(P))\frac{\gamma^{-\frac{1}{p}}}{2\ \mathrm{diam}(F(P))}\varepsilon\\
    &=\frac{1}{2}\varepsilon.
\end{align*}
By the triangle inequality,
\begin{align*}
    \rho_p(f,g)&\le\rho_p(f,P)+\rho_p(P,g)\\
    &\le \gamma^{\frac{1}{p}}\rho_\infty(f,P)+\rho_p(P,g)\\
    &<\frac{1}{2}\gamma^{\frac{1}{p}}\gamma^{-\frac{1}{p}}\varepsilon+\frac{1}{2}\varepsilon
\\
&=\varepsilon
\end{align*}
holds. By the definition of $g$, we have $g\in \per(\frac{d}{dx})$, and the conclusion follows. 
\end{proof}

From the above discussion, the third main theorem holds.

\begin{thm}\label{qwert}
       The differential operator $\frac{d}{dx}$ is chaotic on $(C^\infty(I),\rho_p)$. 
\end{thm}

\begin{proof}
  The conclusion follows from Proposition~\ref{huyu}, Proposition~\ref{haru} and 
Proposition~\ref{denseorbit}. 
\end{proof}

\section{Proof of the Main Theorems on $[a,b]$}

In this section, we prove that the previous propositions also hold for $C^\infty{([a,b])}$.

\subsection{Conjugacy via the Translation Operator}

\begin{dfn}\label{def:T}
  We define linear operators $T,T'$ by
   
    \[T:(C^{\infty}([0,b-a]),\|\cdot\|_p) \to (C^{\infty}([a,b]),\|\cdot\|_p)  ,\quad(Tf)(x):=f(x-a)\]
 \[T':(C^{\infty}([a,b]),\|\cdot\|_p) \to (C^{\infty}([0,b-a]),\|\cdot\|_p) ,\quad(T'f)(x):=f(x+a).\]
\end{dfn}

We will show that $T$ is a conjugacy.

\begin{lem}\label{japan}
    $T'$ is the inverse map of $T$. In particular, $T$ is bijective. 
\end{lem}

\begin{proof}
 For $f\in C^{\infty}([a,b])$ and $x\in[a,b]$, we have 
 \begin{align*}
     (TT'f)(x)=(Tf)(x+a)=f(x+a-a)=f(x).
 \end{align*}
    Thus $TT'f=f$.  This means $TT'=\mathrm{id}_{C^{\infty}([a,b])}$. Similarly, $T'T=\mathrm{id}_{C^{\infty}([0,b-a])}$ holds.
\end{proof}

\begin{lem}\label{china}
    $T$ and $T'$ are isometries. In particular, $T$ and $T'$ are bounded linear operators.
\end{lem}

\begin{proof}
We only show the proof for $T$, since the proof for $T'$ is similar. First, we show the case $p=\infty$. For $f\in C^{\infty}([0,b-a])$, we have
   \begin{align*}
        \|Tf\|_{\infty}&=\sup_{x\in[a,b]}|f(x-a)|=\sup_{x\in[0,b-a]}|f(x)|=\|f\|_\infty .
   \end{align*}
Next, we show the case $1\le p<\infty$. For $f\in C^{\infty}([0,b-a])$, 
\begin{align*}
     \|Tf\|_p^p&=\int_a^b \left|f(x-a)\right|^pdx\\
&=\int_0^{b-a} \left|f(x)\right|^pdx\\
     &=\|f\|_p^p.
\end{align*}
holds. 
\end{proof}

\begin{prop}\label{america}
    The linear operator $T$ is a conjugacy between \\$(C^{\infty}([0,b-a]),\|\cdot\|_p,\frac{d}{dx})$ and $(C^{\infty}([a,b]),\|\cdot\|_p,\frac{d}{dx})$. 
\end{prop}

\begin{proof}
    A bounded linear operator is continuous. Thus, by Lemma~\ref{japan} and Lemma~\ref{china}, $T$ is a homeomorphism. For $f\in C^{\infty}([a,b])$, we have
    \begin{align*}
        \left(\frac{d}{dx}Tf\right)(x)&=\frac{d}{dx}\left(f(x-a)\right)\\
        &=(x-a)'\frac{df}{dx}(x-a)\\
        &=\frac{df}{dx}(x-a)\\
         &=\left(T\frac{df}{dx}\right)(x)\\
         &=\left(T\frac{d}{dx}f\right)(x).\\
    \end{align*}
    Therefore $\frac{d}{dx}T=T\frac{d}{dx}$ holds. 
\end{proof}

\subsection{Main Theorems on $[a,b]$}

By Proposition~\ref{america}, the statements that hold for $C^\infty(I)$ are also valid for $C^\infty([a,b])$. In this subsection, we restate them for $C^\infty([a,b])$. The following is the first main theorem.

\begin{thm}\label{tekoteko}
    There exists an $A\subset C^\infty([a,b])$ such that $(\Lambda,d_{\Lambda},\sigma)$ and $(A,\rho_p,\frac{d}{dx})$ are conjugate.
\end{thm}
\begin{proof}
    The conclusion follows from Theorem~\ref{potato} and Proposition~\ref{america}. 
\end{proof}

\begin{prop}\label{tyotyotyo}
    Let $F\subset\R$ be a finite set with $2\le\#F<\infty$. The set $T(E_F)$ is a $\frac{d}{dx}$-invariant set and $\frac{d}{dx}$ is chaotic on $T(E_F)$. 
\end{prop}

\begin{proof}
     The conclusion follows from Proposition~\ref{kyouyaku} and Proposition~\ref{america}. 
\end{proof}

The following is the second main theorem. 

\begin{thm}\label{kuromichan}
        For any $f\in C^{\infty}([a,b])$ and $\varepsilon>0$, there exists a finite set  $F\subset\R$ with $2\le \#F <\infty$ such that $\rho_p(f,T(E_F))<\varepsilon$. 
\end{thm}

\begin{proof}
    Let $f\in C^{\infty}([a,b])$ and let $\varepsilon>0$. By 
    Theorem~\ref{syut}, we can find a finite set $F\subset \R$ with $2\le\#F<\infty$ such that $\rho_p(T'f,E)<\varepsilon$. By Lemma~\ref{china}, we have
    \begin{align*}
       \rho_p(f,T(E))
        &=\inf_{g\in E}\|f-Tg\|_p\\
        &=\inf_{g\in E}\|T'f-T'Tg\|_p\\
        &=\inf_{g\in E}\|T'f-g\|_p\\
        &=\rho_p(T'f,E)\\
        &<\varepsilon,
    \end{align*}
   and the conclusion follows.
\end{proof}

\begin{rem}
    Theorem~\ref{kuromichan} implies that the following statement holds. \par
      We define $\mathcal{F}=\{F\subset \R\mid 2\le\#F<\infty\}$. Then, $\bigcup_{F\in\mathcal{F}}T(E_F)$ is dense in $C^\infty([a,b])$.
\end{rem}

\begin{cor}\label{kitexichan}
    For any $f\in C^{\infty}([a,b])$, $f$ can be approximated by a filtration $\{T(E_{F(n)})\}_{n=1}^\infty$. That is, there exists a family $\{F(n)\}_{n=1}^{\infty}\subset\R$ satisfying the following conditions.
    \begin{enumerate}
        \item $2\le\#F(n)<\infty$ for any $n\ge 1$
        \item $T(E_{F(n)})\subset T(E_{F(n+1)})$ for any $n\ge1$.
        \item $f\in\overline{\bigcup_{n=1}^\infty T(E_{F(n)})}$, where the closure is with respect to $\rho_p$.
    \end{enumerate}
\end{cor}

\begin{proof}
   The conclusion follows from Corollary~\ref{shimauma} and Proposition~\ref{america}. 
\end{proof}

\begin{rem}
    If $X$ is a space such that $(C^\infty([a,b]),\rho_p)$ can be densely embedded in $X$, a statement similar to Corollary~\ref{kitexichan} holds. That is, for any $f\in X$, there exists a family $\{F(n)\}_{n=1}^{\infty}\subset\R$ satisfying the following conditions.
    \begin{enumerate}
        \item $2\le\#F(n)<\infty$ for any $n\ge 1$
        \item $T(E_{F(n)})\subset T(E_{F(n+1)})$ for any $n\ge1$.
        \item $f\in\overline{\bigcup_{n=1}^\infty T(E_{F(n)})}$, where the closure is with respect to $\rho_p$.
    \end{enumerate}
\end{rem}

The following is the third main theorem.

\begin{thm}\label{kaochan}
      The differential operator $\frac{d}{dx}$ is chaotic on $(C^\infty([a,b]),\rho_p)$. 
\end{thm}

\begin{proof}
    The conclusion follows from Theorem~\ref{qwert} and Proposition~\ref{america}. 
\end{proof}

\subsection{Sensitivity to Initial Conditions}

\begin{dfn}\label{qwkl}
    A discrete dynamical system $(X,f)$ is said to have sensitivity to initial conditions if it satisfies the following: 
    There exists $\beta>0$ such that for any $x\in X$ and any neighborhood $U$ of $x$, there exists $y\in U$ and $n\in \N$ such that $d(f^n(x), f^n(y))>\beta$.
\end{dfn}

The constant $\beta$ in Definition~\ref{qwkl} measures the spread of the error. A natural question arises regarding how large $\beta$ can be for $(C^\infty([a,b]), \frac{d}{dx})$. In fact, $\beta$ can be taken to be arbitrarily large.

\begin{thm}\label{oioi}
      Let $\beta>0$. For any $f\in C^\infty([a,b])$ and any neighborhood $U$ of $f$, there exist $g\in U$ and $n\in \N$ such that $\rho_p(\frac{d^n}{dx^n}f, \frac{d^n}{dx^n}g)>\beta$.
\end{thm}

\begin{proof}
    By Lemma~\ref{china}, Proposition~\ref{america} and Lemma~\ref{pqineq}, it suffices to prove the assertion in $(C^\infty([0,\gamma]),\rho_\infty)$.\par
    We fix $\beta>0$ and $f\in C^\infty([0,\gamma])$. Let $U$ be a neighborhood of $f$. Then, there exists $\varepsilon>0$ such that $B(f,\varepsilon)\subset U$, where $B(f,\varepsilon)$ is the open $\varepsilon$-ball with respect to $\rho_\infty$. We consider two cases. \par
    \textbf{Case 1.} $\sup\{\|\frac{d^k}{dx^k}f\|_\infty\}_{k=0}^\infty=\infty$.\par
    From Lemma~\ref{waieru}, we can find a polynomial function $P$ such that 
    \[\rho_\infty(f,P)<\varepsilon.\]
    This means $P\in B(f,\varepsilon)\subset U$. By the assumption of Case 1, there exists $n\ge\deg(P)+1$ such that $\|\frac{d^n}{dx^n}f\|_\infty>\beta$. Since $\frac{d^n}{dx^n}P=0$, we have
    \begin{align*}
        \rho_\infty\left(\frac{d^n}{dx^n}f, \frac{d^n}{dx^n}P\right)=\rho_\infty\left(\frac{d^n}{dx^n}f,0\right)=\left\|\frac{d^n}{dx^n}f\right\|_\infty>\beta.
    \end{align*}
    Therefore, it suffices to choose $g = P$.\par
    \textbf{Case 2.} $\sup\{\|\frac{d^k}{dx^k}f\|_\infty\}_{k=0}^\infty<\infty$.\par
    We set $M=\sup\{\|\frac{d^k}{dx^k}f\|_\infty\}_{k=0}^\infty$. 
    From Lemma~\ref{waieru}, we can find a polynomial function $P$ such that 
    \[\rho_\infty(f,P)<\frac{\varepsilon}{2}.\]
    We write $P$ as
        \[P(x)=\frac{p_N}{N!}x^N+\frac{p_{N-1}}{(N-1)!}x^{N-1}+\cdots+\frac{p_1}{1!}x+\frac{p_0}{0!},\]
    where $N=\deg(P)$. We set $c=\max\{0,p_0,p_1,\dots,p_N\}+M+\beta$. We define the finite set $F\subset \R$ by 
        \[F=F(P)\cup \{c\}.\]
    Then, it is clear that $2\le \#F<\infty$. By Lemma~\ref{etaxi}, there exists $N'\in\N$ such that 
    \[\zeta_{N'+1}<\frac{\varepsilon}{2~\mathrm{diam}(F)}.\]
   We set $L=\max\{N,N'\}$ and $n=L+1$. For each $j=0,1,2,\dots$, we define $b_j$ by
\begin{align*}
     b_j&:=\left\{ \,
    \begin{aligned}
    & p_j\ &(0\le j \le N) \\
    &  0   &    (j\in \{k\in\N\mid N< k \le L\})  \\
    & c\ &(j> L)  \\
    \end{aligned}
\right.
\end{align*}
and we define $g\in C^\infty(I)$ by $g(x)=\sum_{k=0}^\infty\frac{b_k}{k!}x^k$. Since $P,g\in E_F$, from Lemma~\ref{piano} and Lemma~\ref{etaxi}, we have
\begin{align*}
    \rho_\infty(P,g)\le \mathrm{diam}(F)\zeta_{L+1}\le\mathrm{diam}(F)\zeta_{N'+1}<\frac{\varepsilon}{2}.
\end{align*}
Thus
\begin{align*}
    \rho_\infty(f,g)&\le\rho_\infty(f,P)+\rho_\infty(P,g)\\
    &<\frac{\varepsilon}{2}+\frac{\varepsilon}{2}\\
    &=\varepsilon
\end{align*}
holds. This inequality implies $g\in B(f,\varepsilon)\subset U$. By the definition of $g$, we have
    \begin{align*}
        \left\| \frac{d^n}{dx^n}g\right\|_\infty &=\sup_{x\in[0,\gamma]}\sum_{k=0}^{\infty}\frac{c}{k!}x^k\\
        &>\sup_{x\in[0,\gamma]}c\\
        &=c.
    \end{align*}
By the assumption of Case 2, we have
    \begin{align*}
        \rho_\infty\left(\frac{d^n}{dx^n}f,\frac{d^n}{dx^n}g\right)&=\left\|\frac{d^n}{dx^n}g-\frac{d^n}{dx^n}f\right\|_\infty\\
        &\ge\left\|\frac{d^n}{dx^n}g\right\|_\infty-\left\|\frac{d^n}{dx^n}f\right\|_\infty\\
        &>c-M\\
        &=\max\{0,p_0,p_1,\dots,p_N\}+M+\beta-M\\
        &\ge\beta,
    \end{align*}
and the conclusion follows.
\end{proof}

\begin{cor}\label{kjkj}
    The discrete dynamical system $(C^\infty([a,b]),\frac{d}{dx})$ has sensitivity to initial conditions.
\end{cor}

Originally, Devaney defined chaos in \cite{devaney1989introduction} by transitivity, a density of periodic points, and sensitivity to initial conditions. However, it was shown in \cite{MR1157223} that sensitivity to initial conditions follows from transitivity and the density of periodic points. Therefore, Corollary~\ref{kjkj} follows from Theorem~\ref{oioi}; at the same time, it also follows from Theorem~\ref{kaochan}.

This theorem suggests that the error incurred by repeated differentiation can be significant, even if the approximation error is arbitrarily small in the $L^p$-norm. 

\section{Open problems}

In this section, we state a few questions and open problems that are suggested by the results
of this paper.
\begin{enumerate}
    \item It is natural to investigate the behavior of differential operators on $C^\infty(\Omega)$, where $\Omega$ is an open interval, $\R$, $\mathbb{S}^1$, and so on.
    \item It is also natural to investigate the behavior of partial differential operators.
    \item Do similar statements hold for operators other than $\frac{d}{dx}$? For example, what is known about the chaotic behavior of operators of the form $P(\frac{d}{dx})$, where $P$ is a polynomial?
\end{enumerate}

\bibliographystyle{abbrv}

\end{document}